\documentclass{amsart}
\usepackage{amsmath, amssymb,amsthm,latexsym}
\usepackage{enumerate}

\newtheorem{theorem}{Theorem}[section]

\newtheorem*{main theorem}{Main Result}

\theoremstyle{definition}

\newtheorem{question}[theorem]{Question}
\newtheorem{remark}[theorem]{Remark}

\title{Finite semigroups embed in finitely presented congruence-free monoids}

\author{Victor~Maltcev}

\address{Department of Mathematics and Statistics, Sultan Qaboos University, Al-Khodh 123, Muscat, 
Sultanate of OMAN}

\email{\texttt{victor.maltcev@gmail.com}}


\begin{document}

\maketitle

\begin{abstract}
We prove that every finite semigroup embeds in a finitely presented congruence-free monoid, and pose some questions around the Boone-Higman Conjecture.
\end{abstract}

\section{Introduction}

It is a classical result of Rabin that every countable group embeds in a finitely generated simple group, see~\cite{LS}. Boone and Higman~\cite{Higman} proved that a finitely generated group has soluble word problem if and only if it can be embedded in a simple subgroup of a finitely presented group. This prompted them to raise a question which is now referred to as the Boone-Higman Conjecture:
\begin{question}
Does every finitely generated group with soluble word problem embed in a finitely presented simple group?
\end{question}
This question is still an open problem. Note that the condition of having soluble word problem is crucial, as every simple group has soluble word problem.

The aim of this note is to start understanding the Boone-Higman Conjecture for semigroups, where for the counterpart of simple groups we naturally take congruence-free semigroups. The main result of the note is that every finite semigroup embeds a finitely presented congruence-free monoid. We prove this in Section~\ref{sec:main}, but prior to that let us briefly discuss known embedding theorems in Semigroup Thoery and some questions related to the Boone-Higman Conjecture.

There is the classical Bruck-Reilly extensions~\cite{Howie} by which one can embed a finitely generated (resp. finitely presented) semigroup in a simple finitely generated (resp. finitely presented) semigroup. Byleen also proved that every countable semigroup without idempotents embeds in a $2$-generated simple semigroup without idempotents~\cite{Byleen2}. Furthermore, he proved that every countable semigroup embeds in a $2$-generated bisimple monoid~\cite{Byleen1}. But the corresponding question for finitely presented bisimple semigroups at present seems quite hard:
\begin{question}\label{que:1}
Does every finitely presented semigroup embed in a finitely presented bisimple monoid?
\end{question}
A weaker version of this question is
\begin{question}\label{que:2}
Does every finitely presented semigroup embed in a finitely presented regular monoid?
\end{question}
The author believes there should be some connection between these two questions and the Boone-Higman Conjecture.
In effect Questions~\ref{que:1} and~\ref{que:2} ask whether bisimplicity and regularity are Markov properties (the reader may consult~\cite{Alan} for more information about Markov properties).

The most important result closely related to the Boone-Higman Conjecture is the theorem of Byleen~\cite{Byleen3} that every countable semigroup embeds in a $2$-generated congruence-free semigroup, thus showing that the analogue of Rabin's theorem holds for semigroups.

\section{Main Result}\label{sec:main}

\begin{theorem}
Every finite semigroup embeds in a finitely presented
congruence-free monoid.
\end{theorem}

\begin{remark}
Surprisingly, the required embedding requires some work to do, unlike the case with groups -- every finite group embeds in a corresponding alternating group. Also, our embedding resembles the one from~\cite{BM}, and may be could be used for some other purposes.
\end{remark}

\begin{proof}
Let $S=\{s_1,\ldots,s_n\}$ be a finite semigroup and
$\pi:\{1,\ldots,n\}\times\{1,\ldots,n\}\to\{1,\ldots,n\}$ be the
function defined by $s_is_j=s_{\pi(i,j)}$. The needed finitely
presented congruence-free monoid $M$ to contain $S$, is going to
be constructed as follows: First we introduce new auxiliary
letters $x_1,\ldots,x_{n+1},y_1,\ldots,y_{n+1}$ and then construct
a certain function
$f:\{1,\ldots,n+1\}\times\{1,\ldots,n\}\times\{1,\ldots,n+1\}\to\{0,1\}$
so that our monoid will have a presentation given by the finite
complete system
\begin{align*}
s_is_j &\to s_{\pi(i,j)} & 1\leq i,j\leq n\\
x_is_jy_k &\to f(i,j,k) & 1\leq i,k\leq n+1,\quad 1\leq j\leq n\\
x_iy_j &\to 0 & 1\leq i,j\leq n.
\end{align*}
So that this presentation indeed gives rise to a congruence-free
monoid, we are going to construct $f$ such that the following six
conditions will hold:
\begin{itemize}
\item For every pair
$(i,j)\in\{1,\ldots,n+1\}\times\{1,\ldots,n\}$ there exists $k\leq
n+1$ such that $f(i,j,k)=1$; \item For every pair
$(i,j)\in\{1,\ldots,n\}\times\{1,\ldots,n+1\}$ there exists $k\leq
n+1$ such that $f(k,i,j)=1$; \item For every pair
$(i,j)\in\{1,\ldots,n+1\}\times\{1,\ldots,n\}$ there exists $k\leq
n+1$ such that $f(i,j,k)=0$; \item For every pair
$(i,j)\in\{1,\ldots,n\}\times\{1,\ldots,n+1\}$ there exists $k\leq
n+1$ such that $f(k,i,j)=0$; \item For every two distinct pairs
$(i,j)$ and $(p,q)$ from $\{1,\ldots,n+1\}\times\{1,\ldots,n\}$
there exists $k\leq n+1$ such that $f(i,j,k)\neq f(p,q,k)$; \item
For every two distinct pairs $(i,j)$ and $(p,q)$ from
$\{1,\ldots,n\}\times\{1,\ldots,n+1\}$ there exists $k\leq n+1$
such that $f(k,i,j)\neq f(k,p,q)$.
\end{itemize}
We are left to ensure ourselves that such $f$ indeed exists and
that $f$ with such properties yields a congruence-free monoid. By
priority reasons, we start with the latter. Before that, we remark
that if such $f$ exists, then $M$ is $0$-simple.

\vspace{\baselineskip}

\noindent\textbf{$f$ gives rise to a congruence-free monoid}

\vspace{\baselineskip}

For two distinct elements $u$ and $v$ from $M$, we are going to
prove by induction on $|u|+|v|$ that if $\rho$ is a congruence on
$M$ and $u\rho v$, then $\rho=M\times M$. The base case of
induction -- $|u|+|v|=1$ is obvious.

Now we do the induction step $(<|u|+|v|)\mapsto(|u|+|v|)$. So, let
$u$ and $v$ be two distinct elements of $M$ in their normal forms
with respect to the above finite complete system such that $u\rho
v$. Obviously we may assume that both $u$ and $v$ are non-zero.
Consider first the case when both $u$ and $v$ contain $x_i$'s.
This splits essentially into the following three cases:
\begin{itemize}
\item $u\equiv Ux_is_j$ and $v\equiv Vx_ps_q$. Then there exists
$k\leq n+1$ such that one of $x_is_jy_k$ and $x_ps_qy_k$ is $1$
and the other is $0$, so we may apply induction then. \item
$u\equiv Ux_is_j$ and $v\equiv Vx_p$. Then there exists $k\leq
n+1$ such that $x_iy_jy_k=1$ and then $U=Ux_is_jy_k\rho Vx_py_k=0$
and we may apply induction. \item $u\equiv Ux_i$ and $v\equiv
Vx_p$. If $i=p$, then there exists $k\leq n+1$ such that
$x_is_iy_k=x_ps_iy_k=1$ and then we have that $U\rho V$. Since
$U\not\equiv V$, then we may apply induction. So we may assume
that $i\neq p$. Then there exists $k\leq n+1$ such that one of
$x_is_iy_k$ and $x_ps_iy_k$ is $1$ and the other is $0$, and then
it remains to apply induction.
\end{itemize}

Thus, now we may assume that at least one of $u$ and $v$ contains
none of $x_i$'s; and that at least one of $u$ and $v$ contains
none of $y_i$'s. Assume that at least one of $u$ and $v$ contains
$x_i$'s or $y_i$'s. Say, let $u$ contain $x_i$ as the last or the
penultimate letter. Then
$v\in\{s_1,\ldots,s_n,y_1,\ldots,y_{n+1}\}^{\ast}$ and there are
only the following two cases:
\begin{itemize}
\item $u\equiv Ux_i$. Then $vy_i\rho Ux_iy_i=0$ and since $M$ is
$0$-simple, $\rho=M\times M$. \item $u\equiv Ux_is_j$. Then there
exists $k\leq n+1$ such that $x_is_jy_k=0$ and so $vy_k\rho 0$ and
again by $0$-simplicity of $M$, $\rho=M\times M$.
\end{itemize}

So, we are left with the case when $u,v\in\{s_1,\ldots,s_n\}$:
$u\equiv s_i$ and $v\equiv s_j$ for $i\neq j$. Since there exists
$k\leq n+1$ with the property that one of $x_is_iy_k$ and
$x_is_jy_k$ is $1$ and the other is $0$, we have that $1\rho 0$
and so $\rho=M\times M$.

\vspace{\baselineskip}

\noindent\textbf{Existence of $f$}

\vspace{\baselineskip}

In effect, $f$ is nothing else but coloring the cells of the
`$sxy$' parallelepiped $n\times(n+1)\times(n+1)$ in two colors $0$
and $1$. Then the last two above conditions we imposed on $f$ read
as:
\begin{enumerate}[(1)]
\item For every two distinct `vertical' columns (i.e. those
orthogonal to the `$sx$'-plane) the corresponding $(n+1)$-tuples
are distinct; \item For every two distinct `horizontal' rows (i.e.
those orthogonal to the `$sy$'-plane) the corresponding
$(n+1)$-tuples are distinct.
\end{enumerate}
Now we are going to color the cells corresponding to the slice of
the parallelepiped going through the $k$'th cell on the $s$-axes
and orthogonal to $s$-axes (for $1\leq k\leq n$), this will define
$f$ completely. The $k$'th slice corresponds to a
$(n+1)\times(n+1)$-matrix. Put in the first column of this matrix
the vector $(\underbrace{1,\ldots,1}_{\text{$k$
times}},\underbrace{0,\ldots,0}_{\text{$n+1-k$
times}})^{\mathrm{tr}}$, in the second column the cyclic shift of
the first column by the permutation
$\sigma:(x_1,\ldots,x_{n+1})\mapsto(x_{n+1},x_1,\ldots,x_n)$; in
the third column the shift of the second column by $\sigma$ and so
on. Thus, say, the third slice looks like
\begin{equation*}
\begin{pmatrix}
1 & 0 & 0 & \cdots & 0 & 1 & 1\\
1 & 1 & 0 & \cdots & 0 & 0 & 1\\
1 & 1 & 1 & \cdots & 0 & 0 & 0\\
0 & 1 & 1 & \cdots & 0 & 0 & 0\\
0 & 0 & 1 & \cdots & 0 & 0 & 0\\
\vdots & \vdots & \vdots & \ddots & \vdots & \vdots & \vdots\\
0 & 0 & 0 & \cdots & 1 & 0 & 0\\
0 & 0 & 0 & \cdots & 1 & 1 & 0\\
0 & 0 & 0 & \cdots & 1 & 1 & 1
\end{pmatrix}
\end{equation*}
One easily sees that such a coloring indeed has the required six
conditions.
\end{proof}

\end{document}